\documentclass{amsart}
\usepackage{comment}
\usepackage{MnSymbol}
\usepackage{extarrows}
\usepackage{lscape}
\usepackage[all,cmtip]{xy}

\DeclareMathOperator{\GL}{GL}

\DeclareMathOperator{\Gal}{Gal}

\DeclareMathOperator{\ord}{ord}

\newcommand{\Q}{{\mathbb Q}}
\newcommand{\Z}{{\mathbb Z}}

\newcommand{\F}{{\mathbb F}}

\newcommand{\OO}{{\mathcal O}}

\newcommand{\fp}{\mathfrak{p}}
\newcommand{\fq}{\mathfrak{q}}
\newcommand{\mP}{\mathfrak{P}}

\begin {document}

\newtheorem{thm}{Theorem}

\theoremstyle{definition}

\theoremstyle{remark}

\title[]{On asymptotic Fermat over the $\Z_2$-extension of~$\Q$}

\author{Nuno Freitas}
\address{Departament de Matem\`atiques i Inform\`atica,
Universitat de Barcelona (UB),
Gran Via de les Corts Catalanes 585,
08007 Barcelona, Spain}
\email{nunobfreitas@gmail.com}

\author{Alain Kraus}
\address{Sorbonne Universit\'e,
Institut de Math\'ematiques de Jussieu - Paris Rive Gauche,
UMR 7586 CNRS - Paris Diderot,
4 Place Jussieu, 75005 Paris,
France}
\email{alain.kraus@imj-prg.fr}

\author{Samir Siksek}

\address{Mathematics Institute\\
        University of Warwick\\
        CV4 7AL \\
        United Kingdom}
\email{s.siksek@warwick.ac.uk}

\date{\today}
\thanks{Freitas is supported by a Ram\'on y Cajal fellowship (RYC-2017-22262).
Siksek is supported by
EPSRC grant \emph{Moduli of Elliptic curves and Classical Diophantine Problems}
(EP/S031537/1).}

\begin{abstract}
In a recent work the authors 
prove the effective asymptotic Fermat's Last Theorem for the infinite family of
fields $\Q(\zeta_{2^{r+2}})^+$ where $r \ge 0$.
A crucial step in their proof is the following conjecture of Kraus.
Let $K$ be a number field having odd narrow class number and 
a unique prime $\lambda$ above~$2$. 
Then there are no elliptic curves 
defined over $K$ with conductor~$\lambda$ and a $K$-rational point
of order $2$.
In this note we give a new elementary proof of Kraus' conjecture
that makes use only of basic facts about elliptic curves,
Tate curves and Tate modules.

\bigskip

\noindent {\sc{R\'esum\'e:}} Les auteurs ont  d\'emontr\'e r\'ecemment le th\'eor\`eme de Fermat
asymptotique pour la famille infinie  de corps $\Q\left(\zeta_{2^{r+2}}\right)^+$ avec $r\geq 0$.
Un argument essentiel de la   d\'emonstration est reli\'e \`a  la conjecture suivante de Kraus. Soit $K$ un corps
de nombres ayant un nombre de classes restreint impair et un unique id\'eal premier $\lambda$
au-dessus de $2$. Alors il n'existe pas de courbes elliptiques d\'efinies sur $K$,  de conducteur $\lambda$,
ayant un point d'ordre $2$ rationnel sur $K$. On pr\'esente dans cette note une nouvelle preuve \' el\'ementaire
de la conjecture de Kraus,  en utilisant  seulement des r\'esultats de base sur les courbes elliptiques, qui  concernent
les courbes de Tate et les modules de Tate.
\end{abstract}

\maketitle


\section{Introduction}
Let $K$ be a totally real field, and let $\OO_K$
be its ring of integers. 
The Fermat equation with exponent $p$ over $K$ is the equation
\begin{equation}\label{eqn:Fermat}
a^p+b^p+c^p=0, \qquad a,b,c\in \OO_K.
\end{equation}
A solution $(a,b,c)$ of \eqref{eqn:Fermat} is called trivial
if $abc=0$, otherwise non-trivial.
The \emph{asymptotic Fermat's Last Theorem over $K$}
is the statement that there is a bound $B_K$, depending only on the field $K$, 
such that for all primes $p>B_K$, all solutions to~\eqref{eqn:Fermat} are trivial.
If $B_K$ is effectively computable,
we shall refer to this as the
\emph{effective asymptotic Fermat's Last Theorem over $K$}.
In \cite{FKS} the following two theorems are established.
\begin{thm}\label{thm:main}
Let $K$ be a totally real field satisfying the following two hypotheses:
\begin{enumerate}
\item[(a)] $2$ totally ramifies in $K$;
\item[(b)] $K$ has odd narrow class number.
\end{enumerate}
Then the asymptotic Fermat's Last Theorem holds over $K$.
Moreover, if all elliptic curves over $K$ with full $2$-torsion are modular, then
the effective asymptotic Fermat's Last Theorem holds over $K$.
\end{thm}
Let $r \ge 0$, and let $\zeta_{2^{r+2}}$ be a primitive $2^{r+2}$-th root
of unity. Write $\Q_{r,2}=\Q(\zeta_{2^{r+2}})^+$ for the maximal real
subfield of the cyclotomic field $\Q(\zeta_{2^{r+2}})$.
This is the $r$-th layer of the cyclotomic $\Z_2$-extension
of $\Q$.
\begin{thm}\label{thm:cyclotomic}
The effective asymptotic Fermat's Last Theorem holds
over $\Q_{r,2}$.
\end{thm}
Observe that $\Q_{0,2}=\Q$ and $\Q_{1,2}=\Q(\sqrt{2})$.
Thus Theorem~\ref{thm:cyclotomic} generalizes, albeit asymptotically,
both Fermat's Last Theorem over $\Q$
due to Wiles \cite{Wiles}, and the corresponding theorem over $\Q(\sqrt{2})$
due to Jarvis and Meekin \cite{JarvisMeekin}.
\begin{proof}[Proof of Theorem~\ref{thm:cyclotomic}]
Theorem~\ref{thm:cyclotomic} follows from Theorem~\ref{thm:main}
and the fact that $\Q_{r,2}$ has odd
narrow class number, as shown by Iwasawa \cite{Iwasawa}.
The effectivity follows as elliptic curves
over $\Z_p$-extensions of $\Q$ are modular thanks 
to the work of Thorne~\cite{Thorne}.
\end{proof}
The proof of Theorem~\ref{thm:main}
builds on many deep results, including modularity
lifting theorems over totally real fields due to Kisin, Gee and others,
Merel's uniform boundedness theorem, and Faltings' theorem on rational points
on curves of genus~$\ge~2$, and of course the strategy
of Frey, Serre, Ribet, Wiles and Taylor exploited in Wiles'
proof of Fermat's Last Theorem. A crucial ingredient
in the proof of Theorem~\ref{thm:main}
is furnished by the following theorem,
which had originally been a conjecture of Kraus \cite{Kraus}. 
\begin{thm}\label{thm:Krausgen}
Let $\ell$ be a rational prime. Let $K$ be a number field satisfying
the following conditions:
\begin{enumerate}
\item[(i)] $\Q(\zeta_\ell) \subseteq K$, where $\zeta_\ell$ is a primitive
$\ell$-th root of unity;
\item[(ii)] $K$ has a unique prime $\lambda$ above $\ell$;
\item[(iii)] $\gcd(h_K^+,\ell(\ell-1))=1$ where $h_K^+$ is the 
narrow class number of $K$.
\end{enumerate}
Then there is no elliptic curve $E/K$ with good reduction away
from $\lambda$, potentially multiplicative reduction at $\lambda$,
and a $K$-rational $\ell$-isogeny.
\end{thm}
In the proof of Fermat's Last Theorem, Ribet's Level Lowering
Theorem asserts that the mod $p$ representation of the Frey elliptic curve
arises from a newform of weight $2$ and level $2$. The fact that there
are no such newforms is a seemingly trivial but indeed crucial
step in the proof of Fermat's Last Theorem. 
In the proof of Theorem~\ref{thm:main},
Theorem~\ref{thm:Krausgen}
(with $\ell=2$) plays a similar r\^ole to the absence of newforms
of weight $2$ and level $2$.
For the deduction of Theorem~\ref{thm:main}
from Theorem~\ref{thm:Krausgen} we refer to \cite{FKS}.
The proof of Theorem~\ref{thm:main} in \cite{FKS} makes
heavy use of the theory of $p$-groups and $p$-extensions.
In the present note we give a simpler 
proof of Theorem~\ref{thm:Krausgen},
which uses nothing beyond 
basic facts about elliptic curves,
Tate curves and Tate modules.

\section{Proof of Theorem~\ref{thm:Krausgen}}
Suppose $K$ satisfies conditions (i)--(iii).
In particular there is a unique prime $\lambda$
of $K$ above $\ell$. Let $E/K$
be an elliptic curve with good reduction away from $\lambda$,
potentially multiplicative reduction at $\lambda$. 
We derive a contradiction by studying 
the mod $\ell$ and the $\ell$-adic representations of $E$ (and 
those of a semistable twist). Write 
$G_K=\Gal(\overline{K}/K)$. Denote a decomposition
and inertia subgroups of $G_K$ corresponding to $\lambda$
by $D_\lambda$ and $I_\lambda$ respectively. 

We first show that $E$ has a quadratic twist $F/K$
with conductor $\lambda$ and a $K$-rational $\ell$-isogeny.
Write $\overline{\rho}_{E,\ell}$ for the mod $\ell$
representation of $E$.
By the theory of the Tate curve 
(c.f.\ \cite[Exercises V.5.11 and V.5.13]{SilvermanII}):
$(\overline{\rho}_{E,\ell} \vert_{D_\lambda})^{\mathrm{ss}} \sim
\tau \cdot \chi_{\ell} \oplus \tau$,
where $\chi_{\ell}$ is the modulo $\ell$ cyclotomic character,
and $\tau$ is a character of $D_\lambda$ which is either
trivial or quadratic.
Moreover, the twist $E \otimes \tau$ is an
elliptic curve defined over $K_\lambda$ having split multiplicative
reduction at $\lambda$.
However, by assumption (i), $\chi_{\ell}$ is trivial on $G_K$.
Hence $(\overline{\rho}_{E,\ell} \vert_{D_\lambda})^{\mathrm{ss}} \sim
\tau \oplus \tau$.

As $E$ has a $K$-rational $\ell$-isogeny, the mod $\ell$
representation is reducible, and we can write
\[
\overline{\rho}_{E,\ell} \sim
\begin{pmatrix}
\theta_1 & * \\
0 & \theta_2
\end{pmatrix}
\]
where $\theta_1$, $\theta_2$ are characters $G_K \rightarrow \F_\ell^*$,
and these must satisfy 
$\theta_1 \vert_{I_\lambda}=\theta_2 \vert_{I_\lambda}=\tau \vert_{I_\lambda}$.
As $\tau$ is a quadratic character we see that $\theta_1^2$ and $\theta_1/\theta_2$
are unramified at $\lambda$. 
Since $E/K$ has good reduction away from $\lambda$, by the criterion
of N\'eron--Ogg--Shafarevich
\cite[Proposition IV.10.3]{SilvermanII}
the characters $\theta_1$ and $\theta_2$
are unramified except possibly at $\lambda$ and the infinite places.
We deduce that $\theta_1^2$ and $\theta_1/\theta_2$ are characters of $G_K$
unramified at the finite places having orders dividing $\ell-1$. 
Assumption (iii) immediately implies
 that $\theta_1=\theta_2$ is a quadratic character
of~$G_K$. We let $F$ be the quadratic twist $F=E\otimes \theta_1$.
Note that locally at~$\lambda$ the curve $F/K$ becomes $E \otimes \tau$ and~$\theta_1$ is unramified away from~$\lambda$, hence $F/K$ has conductor $\lambda$,
and
\begin{equation}\label{eqn:n=1}
\overline{\rho}_{F,\ell} \sim
\begin{pmatrix}
1 & * \\
0 & 1 
\end{pmatrix}.
\end{equation}
Write 
$T_{\ell}(F)$ for the $\ell$-adic
Tate module of $E$, and let 
\[
\rho=\rho_{F,\ell^\infty} : G_K \rightarrow \GL(T_{\ell}(F))
\]
be the representation induced by the action of $G_K$. 

As $F$ has multiplicative reduction at $\lambda$, the theory of the
Tate curve tells us \cite[Exercise V.5.13]{SilvermanII} that there is some choice of basis elements 
$P$, $Q \in T_{\ell}(F)$ such
that 
\begin{equation}\label{eqn:inertiaTate}
\rho \vert_{I_\lambda}=\begin{pmatrix}
\chi_{\ell^\infty} & * \\
0 & 1
\end{pmatrix}
\end{equation}
where 
$\chi_{\ell^\infty}: G_K \rightarrow \Z_\ell^\times$ is the 
$\ell$-adic cyclotomic character. 
Fixing this choice of basis $P$, $Q$,
we will show inductively that, as a representation of $G_K$, we have
\begin{equation}\label{eqn:mod2n}
\rho \equiv \begin{pmatrix}
\chi_{\ell^n} & * \\
0 & 1
\end{pmatrix}
\pmod{\ell^n}
\end{equation}
for all $n \ge 1$, where $\chi_{\ell^n}$ is the mod $\ell^n$ 
cyclotomic character. The case $n=1$ is already established in equation~\eqref{eqn:n=1}.

Now suppose $n \ge 2$ and the result holds for~$n-1$. By the inductive hypothesis,
\[
\rho \equiv \begin{pmatrix}
\chi_{\ell^{n}} +\ell^{n-1} \phi & *\\
\ell^{n-1} \psi & 1+\ell^{n-1} \eta
\end{pmatrix} \pmod{\ell^n}
\]
where $\phi$, $\psi$, $\eta$ are functions $G_K \rightarrow \Z/\ell\Z$.
Let $\sigma_1$, $\sigma_2 \in G_K$. Comparing the expressions modulo $\ell^n$
for $\rho(\sigma_1 \sigma_2)$ with $\rho(\sigma_1) \rho(\sigma_2)$
we obtain
\[
\psi(\sigma_1 \sigma_2) \equiv \psi(\sigma_1)\chi_{\ell^{n}}(\sigma_2)+\psi(\sigma_2) \equiv \psi(\sigma_1)+\psi(\sigma_2) \pmod{\ell};
\]
here we have used the fact that $\chi_{\ell^{n}} \equiv \chi_\ell \pmod{\ell}$
and also the fact that $\chi_\ell=1$ as $\Q(\zeta_\ell) \subseteq K$.
Thus $\psi : G_K \rightarrow \Z/\ell\Z$ is an additive character of $G_K$. 
By \eqref{eqn:inertiaTate},
$\psi$ is unramified at $\lambda$, and at all other finite primes
by N\'eron--Ogg--Shafarevich. 
Since $\psi$ has order dividing $\ell$
assumption (iii) allows us to conclude
that $\psi=0$.

Comparing $\rho(\sigma_1 \sigma_2)$ with $\rho(\sigma_1) \rho(\sigma_2)$ 
once more we obtain
\[
\eta(\sigma_1\sigma_2)=\eta(\sigma_1)+ \eta(\sigma_2) \pmod{\ell},
\]
and deduce, as above, that $\eta=0$.
The fact that the determinant of $\rho$ modulo $\ell^n$ must be $\chi_{\ell^n}$ then implies $\phi = 0$, completing the proof of \eqref{eqn:mod2n}.

To complete the proof it remains to demonstrate a contradiction.
One approach is to observe that \eqref{eqn:mod2n} forces $\rho$
to be reducible and to invoke  
 Serre's Open Image Theorem \cite[Chapter IV]{SerreBook}
for a contradiction, 
because $F$ does not have complex multiplication (it has multiplicative reduction prime~$\lambda$).

There is however a more elementary argument which also yields a contradiction.
Let $\fp$ be any prime of $K$ distinct 
from~$\lambda$.
Let $P_n$ and $Q_n$ be the images of $P$, $Q$ in $F[\ell^n]$.
From \eqref{eqn:mod2n}, we note the following.
\begin{itemize}
\item The cyclic subgroup $\langle P_n \rangle$ is fixed by $G_K$
and therefore the isogenous elliptic curve $F_n=F/\langle P_n \rangle$
is defined over $K$.
\item $\sigma(Q_n)=a_\sigma P_n +Q_n$
for any $\sigma \in G_K$ where
$a_\sigma \in \Z/\ell^n\Z$. Thus $Q_n+ \langle P_n \rangle$
is a $K$-point of order $\ell^n$ on $F_n$.
\end{itemize}
Since $F_n$ has good reduction at~$\fp$,
by the injectivity of torsion under reduction 
we see that $\ell^n \mid \# F_n(\F_\fp)$ and as $F$ and $F_n$
are isogenous, $\ell^n \mid \# F(\F_\fp)$. This gives
a contradiction for $n$ large.

\end{document}